\documentclass[3p,times]{elsarticle}

\usepackage{amsmath}
\usepackage{amssymb}
\usepackage{graphics}
\usepackage{epsfig}
\usepackage{epstopdf}
\usepackage{sidecap}
\usepackage{subcaption}
\usepackage{algorithm}
\usepackage{algpseudocode}
\usepackage{multicol}
\usepackage{enumitem}
\usepackage{tabularx}
\usepackage{comment}
\usepackage{color}
\usepackage[T1]{fontenc}

\usepackage[figuresright]{rotating}

\def\bs{\boldsymbol}

\begin{document}

\begin{frontmatter}




\title{Dimensionality Reduction of Collective Motion by Principal Manifolds}

\author[label1]{Kelum Gajamannage}
\ead{dineshhk@clarkson.edu}

\author[label2]{Sachit Butail}
\ead{sb4304@nyu.edu}
\author[label2]{Maurizio Porfiri}
\ead{mporfiri@poly.edu}
\author[label1]{Erik M. Bollt\corref{cor1}}
\ead{bolltem@clarkson.edu}
\address[label1]{Department of Mathematics, Clarkson University, Potsdam, NY-13699, USA.}
\address[label2]{Department of Mechanical and Aerospace Engineering, New York University Polytechnic School of Engineering, Brooklyn, NY-11201, USA.}

\cortext[cor1]{Corresponding author}
\begin{abstract}
While the existence of low-dimensional embedding manifolds has been shown in patterns of collective motion, the current battery of nonlinear dimensionality reduction methods are not amenable to the analysis of such manifolds.  This is mainly due to the necessary spectral decomposition step, which limits control over the mapping from the original high-dimensional space to the embedding space. Here, we propose an alternative approach that demands a two-dimensional embedding which topologically summarizes the high-dimensional data. In this sense, our approach is closely related to the construction of one-dimensional principal curves that minimize orthogonal error to data points subject to smoothness constraints. Specifically, we construct a two-dimensional principal manifold directly in the high-dimensional space using cubic smoothing splines, and define the embedding coordinates in terms of geodesic distances. Thus, the mapping from the high-dimensional data to the manifold is defined in terms of local coordinates. Through representative examples, we show that compared to existing nonlinear dimensionality reduction methods, the principal manifold retains the original structure even in noisy and sparse datasets. The principal manifold finding algorithm is applied to configurations obtained from a dynamical system of multiple agents simulating a complex maneuver called predator mobbing, and the resulting two-dimensional embedding is compared with that of a well-established nonlinear dimensionality reduction method.
\end{abstract}


\end{frontmatter}


\section{Introduction}
With advancements in data collection and video recording methods, high-volume datasets of animal groups, such as fish schools \cite{partridge1982structure, Gerlai2010b}, bird flocks \cite{Nagy2010, Ballerini2008}, and insect and bacterial swarms \cite{Branson2009, zhang2009swarming}, are now ubiquitous. However, analyzing these datasets is still a nontrivial task, even when individual trajectories of all members are available. A desirable step that may ease the experimenter's task of locating events of interest is to identify coarse observables \cite{Kolpas2007, Miller2008, Frewen2011} and behavioral measures \cite{Miller2011a} as the group navigates through space. In this context, Nonlinear Dimensionality Reduction (NDR) offers a large set of tools to infer properties of such complex multi-agent dynamical systems. 

Traditional Dimensionality Reduction (DR) methods based on linear techniques, such as Principal Components Analysis (PCA), have been shown to possess limited accuracy when input data is nonlinear and complex \cite{van2009dimensionality}. DR entails finding the axes of maximum variability \cite{kirby2000geometric} or retaining the distances between points \cite{cox2010multidimensional}. Multi Dimensional Scaling (MDS) with Euclidean metric is another DR method which attains low-dimensional representation by retaining the pairwise distance of points in low dimensional representations \cite{cox2010multidimensional}. However, Euclidean distance calculates the shortest distance between two points on a manifold instead of the genuine manifold distance, which may lead to difficulty of inferring low-dimensional embeddings. The isometric mapping algorithm (Isomap) resolves the problem associated with MDS by preserving the pairwise geodesic distance between points \cite{tenenbaum2000global}; it has recently been used to analyze group properties in collective behavior, such as the level of coordination and fragmentation \cite{abaid2012topological, delellis2014collective, arroyo2012reverse}.  Within Isomap, however, short-circuiting \cite{samko2006selection} created by faulty connections in the neighborhood graph, manifold non-convexity \cite{tipping2001sparse, DeLellis2013} and holes in the data \cite{li2005supervised} can degrade the faithfulness of the reconstructed embedding manifold. 

Diffusion maps \cite{coifman2006diffusion} have also been shown to successfully identify coarse observables  in collective phenomena \cite{aureli2012portraits} that would otherwise require hit-and-trial guesswork \cite{frewen2011coarse}. Beyond Isomap and diffusion maps, the potential of other NDR methods to study collective behavior is largely untested. For example, Kernel PCA (KPCA) requires the computation of the eigenvectors of the kernel matrix instead of the eigenvectors of the covariance matrix of the data \cite{shawe2004kernel} but this is  computationally expensive \cite{van2009dimensionality}. Local Linear Embedding (LLE) embeds high-dimensional data through global minimization of local linear reconstruction errors \cite{van2009dimensionality}.  Hessian LLE (HLLE) minimizes the curviness of the higher dimensional manifold by assuming that the low-dimensional embedding is locally isometric \cite{donoho2003hessian}. Laplacian Eigenmaps (LE) perform a weighted minimization (instead of global minimization as in LLE) of the distance between each point and its given nearest neighbors to embed high dimensional data \cite{belkin2001laplacian}. 

Iterative NDR approaches have also been recently developed in order to bypass spectral decomposition which is common in most of NDR methods \cite{gashler2007iterative}. Curvilinear Component Analysis (CCA) employs a self-organized neural network to perform two tasks, namely, vector quantization of submanifolds in the input space and nonlinear projection of  quantized vectors onto a low dimensional space \cite{demartines1997curvilinear}. This method minimizes the distance between the input and output spaces. Manifold Sculpting (MS) transforms data by balancing two opposing heuristics: first, scaling information out of unwanted dimensions, and second, preserving local structure in the data. The MS method, which is robust to sampling issues, iteratively reduces the dimensionality by using a cost function that simulates the relationship among points in a local neighborhoods\cite{gashler2007iterative}. The Local Spline Embedding (LSE) is another NDR method which embeds the data points using splines that map each local coordinate into a global coordinate of the underlying manifold by minimizing the reconstruction error of the objective function \cite{xiang2009nonlinear}. This method reduces the dimensionality by solving an eigenvalue problem while the local geometry is exploited by the tangential projection of data.  LSE  assumes that the data is not only unaffected by noise or outliers, but also, sampled from a smooth manifold which ensures the existence of a smooth low dimensional embedding.

Due to the global perspective of all these methods, none of them provide sufficient control over the mapping from the original high-dimensional dataset to the low-dimensional representation, limiting the analysis in the embedding space. In other words, the low-dimensional coordinates are not immediately perceived as useful, whereby one must correlate the axes of the embedding manifold with selected functions of known observables to deduce their physical meaning \cite{frewen2011coarse, tenenbaum2000global}. In this context, a desirable feature of DR that we  emphasize here is the regularity in the spatial structure and range of points on the embedding space, despite the presence of noise.

\begin{figure}[ht!]
	\centering
	\includegraphics[width=.995\linewidth]{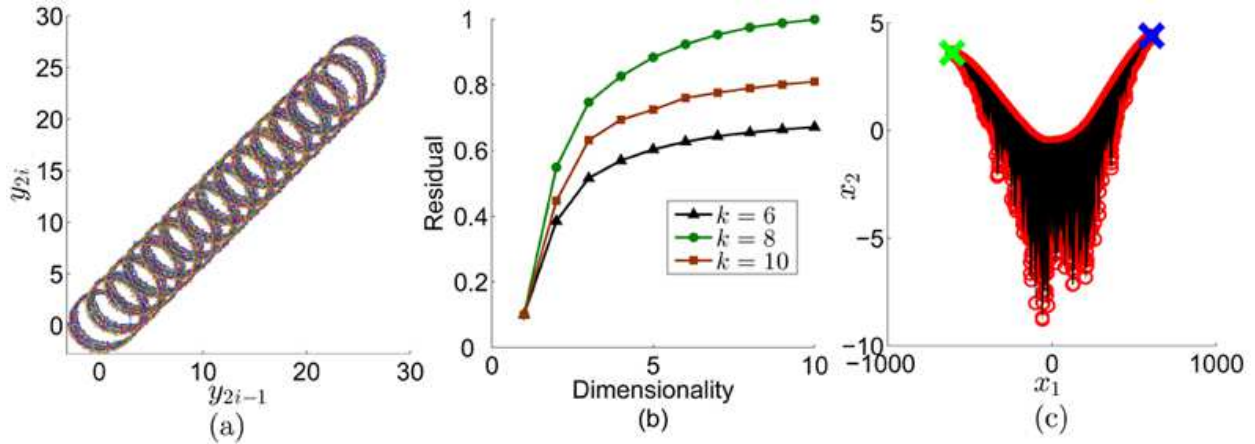}
	\caption{Using Isomap to create a two-dimensional embedding of a simulation of collective behavior. (a) Predator mobbing of twenty agents moving on a translating circular trajectory on a plane (enclosing a predator moving at constant speed at a 45$^{\circ}$ angle), axes $y_{2i-1}$, $y_{2i}$ generally represent coordinates for the $i$-th agent. (b) Scaled residual variance of candidate low-dimensional embeddings produced by Isomap using different nearest neighbor values $k$ (green-circle, brown-square, and black-triangle), and (c) two-dimensional representation of the data for five nearest neighbors (black-triangle). Green and blue crosses mark the start and end points of the trajectory.}
	\label{fig:predator_mobbing_isomap}
\end{figure}

With regard to datasets of collective behavior, nonlinear methods have limited use for detailed analysis at the level of the embedding space. This is primarily because the majority of these methods collapse the data onto a lower dimensional space, whose coordinates are not guaranteed to be linear functions of known system variables \cite{Lee2004a}.  In an idealized simulation of predator induced mobbing \cite{dominey1983mobbing}, a form of collective behavior where a group of animals crowd around a moving predator, two degrees of freedom are obvious, namely, the translation of the group and the rotation about the predator (center of the translating circle). This two-dimensional manifold is not immediately perceived by Isomap, even for the idealized scenario presented in Figure \ref{fig:predator_mobbing_isomap}, where a group of twenty individuals rotate about a predator moving at a constant speed about a line bisecting the first quadrant. Specifically, the algorithm is unable to locate a distinct elbow in the residual variance vs. dimensionality curve, notwithstanding substantial tweaking of the parameters. Thus,  the inferred dimensionality is always 1 (Fig. \ref{fig:predator_mobbing_isomap}b). For a two-dimensional embedding (Fig. \ref{fig:predator_mobbing_isomap}c), visual comparison of the relative scale of the axes indicates that the horizontal axis represents a greater translation than the vertical axis. It is likely that the horizontal axis captures the motion of the group along the translating circle. The vertical axis could instead be associated with  (i) motion about the center of the circle, or (ii) noise, which is independent and identically distributed at each time step. The ambiguity in determining the meaning of such direction indicates a drawback of Isomap in provide meaningful interpretations of the low-dimensional coordinates.

An alternative approach to DR, one that does not require heavy matrix computations or orthogonalization, involves working directly on raw data in the high-dimensional space \cite{Ball1965isodata, Hastie1989}. We propose here a method for DR that relies on geodesic  rather than Euclidean distance and emphasizes manifold regularity. Our approach is based on a spline representation of the data that allows us to control the expected manifold regularity. Typically, this entails conditioning the data so that the lower dimensions are revealed. For example, in \cite{Ball1965isodata}, raw data is successively clustered through a series of lumping and splitting operations until a faithful classification of points is obtained. In \cite{Hastie1989}, one-dimensional parameterized curves are used to summarize the high-dimensional data by using a property called self-consistency, in which points on the final curve coincide with the average of raw data projected on itself. In a similar vein, we construct a PM of the high-dimensional data using cubic smoothing splines.
Summarizing the data using splines to construct a PM of the data has two advantages: (i) it respects the data regularity by enabling access to every point of the dataset, and (ii) it gives direct control over the amount of noise rejection in extracting the true embedding through the smoothing parameter. We illustrate the algorithm using the standard swiss roll dataset, typically used in NDR methods.  Then we validate  the method on three different datasets, including an instance of collective behavior similar to that in Figure \ref{fig:predator_mobbing_isomap}.

This paper is organized as follows. Section \ref{sec:algorithm} describes the three main steps of the algorithm using the swiss roll dataset as an illustrative example. In Section \ref{sec:examples}, we validate and compare the algorithm on a paraboloid, a swiss roll with high additive noise, and a simulation of collective animal motion. Section \ref{sec:performance} analyzes the performance of the algorithm by comparing topologies in the original and embedding space as a function of the smoothing parameter, noise intensity, and data density. We conclude in Section \ref{sec:conclusion} with a discussion of the algorithm performance and ongoing work. Individual steps of the algorithm are detailed in \ref{sec:algorithms}. Computational complexity is discussed in \ref{sec:complexity}.

\section{Dimensionality reduction algorithm to find principal manifold} \label{sec:algorithm}
Dimensionality reduction is defined as a transformation of high-dimensional data into a meaningful representation of reduced dimensionality \cite{van2009dimensionality}. In this context,  a $d$-dimensional input dataset is embedded onto an $e$-dimensional space such that $e<<d$. The mapping from point $\bs{x} \in \mathbb{R}^e$ in the embedding space to the corresponding point $\bs{y} \in \mathbb{R}^d$ in the original dataset is
\begin{equation}
\Phi:\mathbb{R}^e\rightarrow \mathbb{R}^d
\end{equation}

Differently from other approaches, the proposed DR algorithm is based on a PM of the raw dataset, obtained by constructing a series of cubic smoothing splines on slices of data that are partitioned using locally orthogonal hyperplanes. The resulting PM summarizes the data in the form of a two-dimensional embedding, that is, $e=2$. The PM finding algorithm proceeds in three steps: clustering, smoothing, and embedding. The clustering step partitions the data using reference points into non-overlapping segments called slices. This is followed by locating the longest geodesic within each slice. A cubic smoothing spline is then fitted to the longest geodesic. A second set of slices, and corresponding splines, are created in a similar way resulting in a two-dimensional PM surface. The PM is constructed in the form of intersection points between pairs of lines, one from each reference point. Any point in the input space is projected on the PM in terms of the intersection points and their respective distances along the cubic splines. Table \ref{tab:list_of_symbols} lists the notation used in this paper.

\begin{table}
	\centering
	\caption{Nomenclature}
	\begin{tabular}{c | p{10cm}}
		Notation & Description\\
		\hline
		 &  \\
		$\mathcal{D}$ & Input dataset\\
		$d$ & Input dimension\\
		$\bs{y}$ & A point in the input dataset \\
		$n$ & Number of points in the input dataset \\
		$e$ & Output embedding dimension (equal to two)  \\
		$\bs{x}$ & A point in the embedding space \\
		$\bs{z}$ & An external input point for the embedding \\		
		$\bs{\mu}$ & Mean (centroid) of the input dataset \\
		$p$ & Smoothing parameter \\	
		$k$ & Number of neighbors in the nearest neighbor search  \\
		$\bs{O}$ & Origin of the underlying manifold \\	
		$\bs{q}_{1,2}$ & Reference points 1 and 2 \\
		$\mathcal{C}^{1,2}_j$ & $j$-th cluster for reference points 1 and 2\\
		$n_\mathcal{C}^{1, 2}$ & Number of clusters for reference points 1 and 2\\
		$\bs{S}^{1,2}_j$ & $j$-th cubic spline for reference points 1 and 2 \\
		$\mathcal{G}^{1,2}_j$ & $j$-th geodesic for reference points 1 and 2 \\
		$\bs{t}^{l,m}$ & Intersection point in $d$-dimensional space between splines $l$ and $m$ \\							$(\bs{u}^1_{\bs{z}}, \bs{u}^2_{\bs{z}})$ & Tangents to splines at the point $\bs{z}$ \\
		$(\bs{v}_1, \sigma_1), (\bs{v}_2, \sigma_2)$ & First two Principal Components (PCs) of the input dataset, where $\bs{v}_i$ is the $d$-dimensional coefficients and $\sigma_i$ is the eigenvalue of the $i$-th PC \\		
		
	\end{tabular}
	\label{tab:list_of_symbols}

\end{table}

\subsection{Clustering}
\begin{figure}[!ht]
	\centering
	\includegraphics[width=1\textwidth]{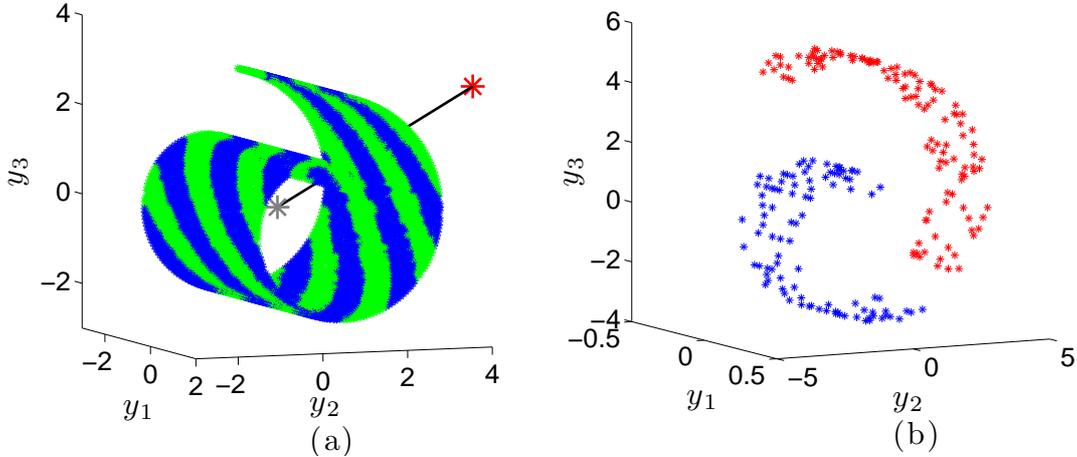}	        	
	\caption{Clustering the raw data into slices. (a) Using a single reference point (red) located far from the data points, the dataset can be sliced into non-overlapping sections which are perpendicular to the line joining the centroid (gray star) and the reference point (red star), then sections are used to draw a grid pattern using splines. (b) A hole in the data will result in multiple sub-clusters (shown in different colors red and blue) within a cluster.}
	\label{fig:sliced_swiss_roll}
\end{figure}

Consider an input dataset $\mathcal{D}$ described by $n$ $d$-dimensional points $\bs{y}_1, \ldots, \bs{y}_n$. Clustering begins by partitioning the data into non-overlapping clusters. The partitioning is performed by creating hyperplanes orthogonal to the straight line joining an external reference point $\bs{q}\in \mathbb{R}^d$ through the mean $\bs{\mu}\in \mathbb{R}^d$ to some point $\tilde{\bs{q}}\in \mathbb{R}^d$. Although any reference point may be selected, we use Principal Components (PCs) to make the clustering procedure efficient in data-representation. In particular, Principal Component Analysis (PCA) of the matrix $\mathcal{D}=\begin{bmatrix} \bs{y}_1 | \bs{y}_2 | \ldots | \bs{y}_n \end{bmatrix}$ is performed to obtain two largest principal components $(\bs{v}_1, \sigma_1)$ and $(\bs{v}_2, \sigma_2)$, where $\bs{v}_i$ is the $d$-dimensional coefficients and $\sigma_i$ is the eigenvalue of the $i$-th PC. The first reference point is chosen such that $\bs{q}_1=\bs{\mu}+\bs{v}_1\sigma_1$. To cover the full range of the dataset, the line joining the points $\bs{q}_1$ and the point $\tilde{\bs{q}_1}=\bs{\mu}- \bs{v}_1\sigma_1$ is divided into $n_{\mathcal{C}}^1$ equal parts using the ratio formula in a $d$-dimensional space \cite{Protter1988}. The $j$-th point $\bs{a}_j \in \mathbb{R}^d$ on this line is 
\begin{equation}
\bs{a}_j=\bs{\mu}+ \frac{\bs{v}_i\sigma_i(2j-n_\mathcal{C}^i)}{n_\mathcal{C}^i}; \ \ j=0, \dots, n_\mathcal{C}^i
\label{eqn:line_partition}
\end{equation}
where $i=1$ for this case of dealing with the first reference point. All points between hyperplanes through points $s_{j-1}$ and $s_j$ are assigned to the cluster $\mathcal{C}_j^1$, $j=1,\ldots,n_\mathcal{C}^1$, by using the following inequality for $k=1, \dots, n$ with $i=1$
\begin{equation}\label{eqn:enq_clustering}
\begin{split}
\cos^{-1}\left(\frac{\left(\bs{y}_k-\bs{a}_{j-1}\right)^\mathrm{T} \left(\bs{q}_i-\bs{a}_{j-1}\right)}{\|\bs{y}_k-\bs{a}_{j-1}\|\|\bs{q}_i-\bs{a}_{j-1}\|}\right)<\pi/2 \hspace{10pt} , \hspace{10pt}
\cos^{-1}\left(\frac{\left(\bs{y}_k-\bs{a}_{j}\right)^\mathrm{T}  \left(\bs{q}_i-\bs{a}_{j}\right)}{\|\bs{y}_k-\bs{a}_{j}\|\|\bs{q}_i-\bs{a}_{j}\|}\right)\ge\pi/2
\end{split}
\end{equation} 
(Fig. \ref{fig:sliced_swiss_roll}a). The clustering method is adapted to detect gaps or holes in the data by setting a threshold on the minimum distance between neighboring points based on a nearest neighbor search \cite{Nene1997simple}. If a set of points do not satisfy the threshold, they are automatically assigned another cluster giving rise to sub-clusters \cite{MacQueen1967} (Fig. \ref{fig:sliced_swiss_roll}b). 

The above procedure is repeated for the second reference point $\bs{q}_2=\bs{\mu}+ \bs{v}_2\sigma_2$ to partition the dataset along the line joining $\bs{\mu}\pm \bs{v}_2\sigma_2$ using the equation (\ref{eqn:line_partition}) with $i=2$. Then, equation (\ref{eqn:enq_clustering}) is used  for $k=1, \dots, n$ with $i=2$ to make another set of clusters $\mathcal{C}_j^2 ; j=1,\ldots,n_\mathcal{C}^2$. The resulting clustering algorithm inputs the data $\mathcal{D}$, and the number of clusters $n_{\mathcal{C}}^{1, 2}$ with respect to each reference point, which can be set on the basis of data density in each direction. The set of clusters for each reference point are stored for subsequent operations. The clustering algorithm is detailed in \ref{sec:algorithms} as algorithm \ref{alg:cluster}.

\subsection{Smoothing} 
The clustering step is followed by smoothing, where cubic splines are used to represent the the clusters. This approach of using a spline to approximate the data is similar to forming a principal curve on a set of points \cite{Hastie1989, Biau2012}. 
Briefly, the principal curve runs smoothly through the middle of the dataset such that any point on the curve coincides with the average of all points that project orthogonally onto this point \cite{Hastie1989}. The algorithm in \cite{Hastie1989} searches for a parameterized function that satisfies this condition of self-consistency. In the context of smoothing splines on a cluster with $m$ points, $\bs{y}_1, \ldots, \bs{y}_{m}$, the principal curve $\bs{g}(\lambda)$ parameterized in $\lambda \in [0,1]$ minimizes over all smooth functions $\bs{g}$ \cite{Hastie1989}
\begin{equation}
\sum_{i=1}^{m}\|\bs{y}_i-\bs{g}(\lambda_i)\|^2 + \kappa\int_0^1\|\bs{g}''(\lambda)\|^2d\lambda,
\label{eqn:pcurve}
\end{equation}
where $\kappa$ weights the smoothing of the spline and $\lambda_i \in [0,1]$ for $i=1, \ldots, m$. Equation (\ref{eqn:pcurve}) yields $d$ individual expressions, one for each dimension in the input space.

\begin{figure}[!ht]
        	\centering
        	\includegraphics[width=1\textwidth]{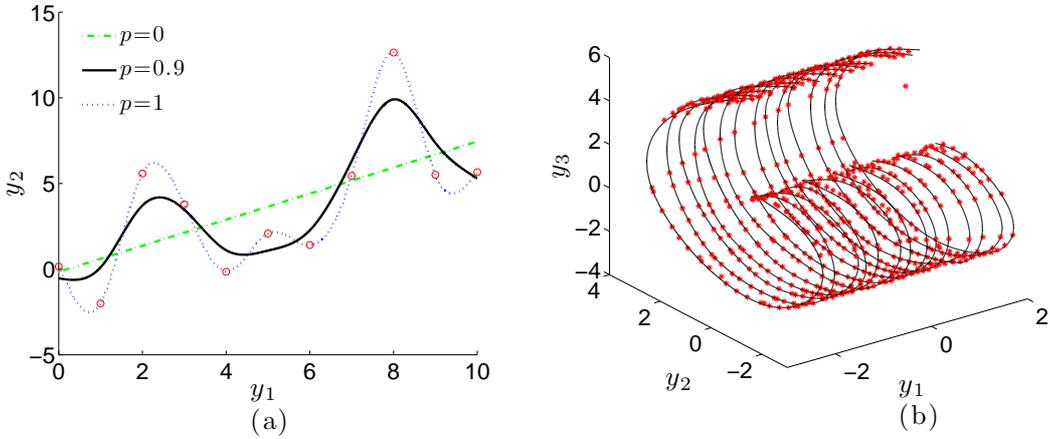}       
           \caption{Approximating clusters with smoothing splines. (a) A cubic smoothing spline is fit onto two-dimensional data using different values of the smoothing parameter $p$. (b) Smoothing splines fit clusters of a swiss roll dataset from a single reference point.}
	\label{fig:p_cubic_smoothing_splines}       	
\end{figure}

Since the spline is created in a multi-dimensional space, the points are not necessarily arranged to follow the general curvature of the dataset. In order to arrange the points that are input to the equation (\ref{eqn:pcurve}), we perform an additional operation to build geodesics within the cluster. The geodesics are created using a range-search\footnote{A neighbor search which choses all neighbors within a specific distance.} with range distance set as the cluster width $2\sigma_1/n_{\mathcal{C}}^1$ or $2\sigma_2/n_{\mathcal{C}}^2$ depending on the reference point. The longest geodesic within the cluster is  used to represent the full range and curvature of the cluster. Using $n_g$ number of $d$-dimensional points $\left\{\bs{y}_i=\left[y_{i,1}, \dots, y_{i,d}\right]^{\mathrm{T}} \vert i=1, \dots, n_g\right\}$ on the longest geodesic, the $d$-dimensional  cubic smoothing spline $\bs{S}(\lambda)=\left[S_1(\lambda) | S_2(\lambda) | \dots | S_d(\lambda)\right]$ parameterized in $\lambda \in [0, 1]$ is computed as a piecewise polynomial that minimizes the quantity \cite{Hastie1989, Bollt2007}
\begin{equation}
	p\sum_{i=1}^{n_g}|y_{i,j}-S_j(\lambda_i)|^2+(1-p)\int_{\lambda_1}^{\lambda_{n_g}}\left[S''_j(\lambda)\right]^2d\lambda,
	\label{eqnCSS}
\end{equation}
where $p \in [0,1]$, in each dimension $j=1, \dots, d$. Thus, $p$ weights the sum of square error and $(1-p)$ weights the smoothness of the spline; conversely $p=0$ gives the natural cubic smoothing spline and $p=1$ gives an exact fit (Fig. \ref{fig:p_cubic_smoothing_splines}a). This procedure is repeated for all clusters for each reference point (Fig. \ref{fig:p_cubic_smoothing_splines}b) resulting in a grid-like surface of smoothing splines. We denote cubic smoothing splines made with respect to the first reference point by $\bs{S}^1_j ; j=1, \dots, n^1_\mathcal{C}$, and second reference point by $\bs{S}^2_j ; j=1, \dots, n^2_\mathcal{C}$. The smoothing algorithm is described in \ref{sec:algorithms} as algorithm \ref{alg:smooth}.

\subsection{Embedding} 
Once constructed, intersection points between all pairs of splines, one from each reference point are located as landmarks to embed new points onto the PM.  Since the splines from two reference points may not actually intersect each other unless $p=1$ for all splines, virtual intersection points between two non-intersecting splines are located at the midpoint on the shortest distance between them. The embedding algorithm loops through all pairs of splines to locate such points, and once located they are denoted as $\bs{t}^{l,m}\in \mathbb{R}^d$ where $(l, m) \in \bs{S}^1_l \times \bs{S}^2_m$.

\begin{figure*}[!ht]
	\includegraphics[width=0.995\linewidth]{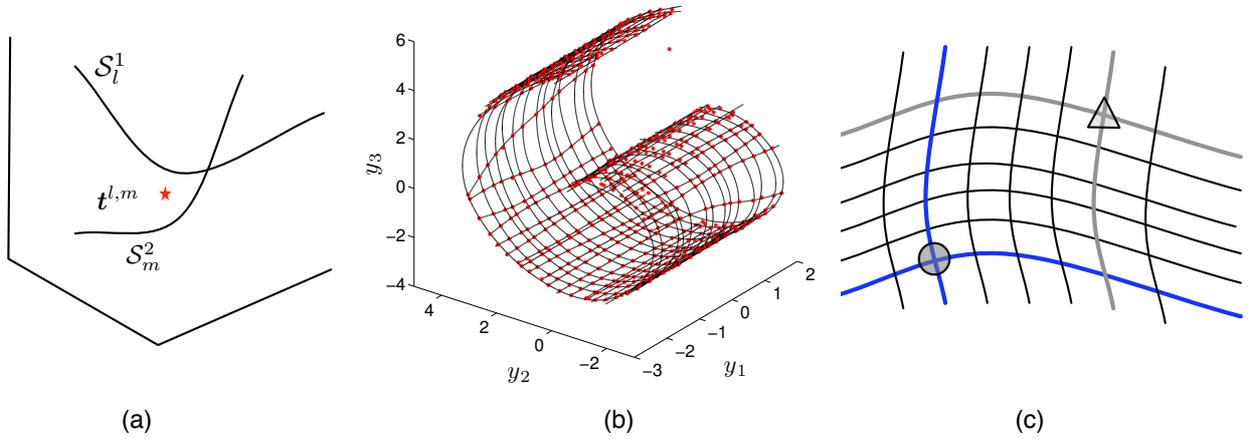}
        \caption{Embedding the data into two-dimensional coordinates. (a) The intersecting point between two splines lies midway on the shortest distance between them (b) intersecting points for a swiss roll and (c) the two-dimensional coordinate of a point (triangle) is the geodesic distance along the intersecting splines (grey) to the axis splines (blue) from an origin (circle).}
         \label{fig:approximating_intersections}
\end{figure*}

In order to define the embedding coordinates, an intersection point is selected randomly and assigned as the manifold origin $O$ of the PM. The set of smoothing splines corresponding to the origin are called axis splines. The embedding coordinates of a point $\bs{z} \in \mathbb{R}^d$ are defined in terms of the spline distance of the intersecting splines of that point from the axis splines (Fig. \ref{fig:approximating_intersections}c). For that, first we find the closest intersection $\bs{\tilde{z}}$ for $\bs{z}$ in Euclidean distance as
\begin{equation}\label{eqn:closest_intersection}
\tilde{\bs{z}} = \mathrm{argmin}_{l,m}\|\bs{t}^{l,m}-\bs{z}\|
\end{equation}
and tangents to splines at this point are called the local spline directions and denoted by $(\bs{u}^1_{\bs{z}}, \bs{u}^2_{\bs{z}})$. Distance from $O$ to $\tilde{\bs{z}}$ is computed along the axis splines by
\begin{equation}\label{eqn:axis_distance}
\left[\tilde{x}_1, \tilde{x}_2\right]^{\mathrm{T}}=\begin{bmatrix}d_l(\tilde{\bs{z}}-O), d_m(\tilde{\bs{z}}-O)\end{bmatrix}^\mathrm{T} 
\end{equation}
where $d_l(b_1-b_2)$ is the distance between points $b_1$ and $b_2$ along the spline $l$. We project the vector $\tilde{\bs{z}}-\bs{z}$ on the local tangential directions $(\bs{u}^1_{\bs{z}}, \bs{u}^2_{\bs{z}})$ as 
\begin{equation}\label{eqn:local_distance}
\begin{split}
\delta x_1= \left(\tilde{\bs{z}}-\bs{z}\right)^\mathrm{T} \bs{u}^1_{\bs{z}}\\
\delta x_2 =\left(\tilde{\bs{z}}-\bs{z}\right)^\mathrm{T} \bs{u}^2_{\bs{z}}
\end{split}
\end{equation}
at the intersection point $\tilde{\bs{z}}$. The final embedding coordinates are obtained by summing the quantities in equations (\ref{eqn:axis_distance}) and (\ref{eqn:local_distance}) as
\begin{equation}
\left[x_1, x_2\right]^{\mathrm{T}}=\left[\tilde{x}_1+ \delta x_1, \tilde{x}_2+ \delta x_2\right]^{\mathrm{T}}
\end{equation} 
The embedding algorithm is detailed in \ref{sec:algorithms} as algorithm \ref{alg:embed}.

\subsection{Mapping between the manifold and raw data}\label{subsec:invMap}
\begin{figure}[htp]
\centering`
	\includegraphics[width=0.4\linewidth]{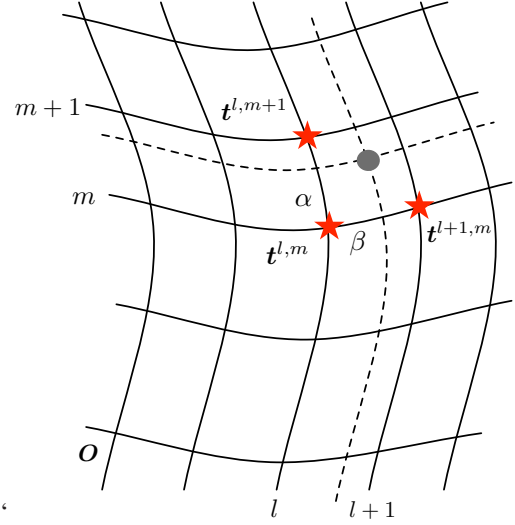}	
        \caption{The mapping from the embedding manifold to the original high-dimensional space is found by finding the projections $\alpha, \beta$ on the local coordinate space and propagating the same into the high dimensions. Local coordinate space can be created using adjacent splines ($l, l+1$), and ($m, m+1$).}
        \label{fig:inverse_mapping_interpolation}
\end{figure} 

The mapping between the manifold and the raw data is created by inverting the steps of the embedding algorithm. In particular, a two-dimensional point $\bs{x}=\begin{bmatrix}x_1, x_2\end{bmatrix}^\mathrm{T}$ in the embedding space is located in terms of the
embedding coordinate $\bs{x}^{l,m}$ of the closest intersection point $\bs{t}^{l,m} \in \mathbb{R}^d$. We then complete the local coordinate system by finding the two adjacent intersection points. The resulting three points on the PM, the origin and the two adjacent points denoted by $\{\bs{x}^{l,m}, \bs{x}^{l,m+1}, \bs{x}^{l+1,m}\}$, are used to solve for ratios $\alpha$ and $\beta$ such that
\begin{equation}
	\begin{split}
		\alpha&=\frac{x_2-x^{l,m}_2}{x_2^{l,m+1}-x^{l,m}_2} \\
		\beta&=\frac{x_1-x_1^{l,m}}{x_1^{l+1,m}-x_1^{l,m}}.
	\end{split}
\end{equation}
The ratios are then propagated to the high-dimensional space (Fig. \ref{fig:inverse_mapping_interpolation}) by using the same relation on the corresponding points in the higher dimension $\{\bs{t}^{l,m}, \bs{t}^{l,m+1}, \bs{t}^{l+1,m}\}$ to obtain the high-dimensional point 
\begin{equation}
		\bs{y}=\bs{t}^{l,m}+(\bs{t}^{l,m+1}-\bs{t}^{l,m})\alpha+(\bs{t}^{l+1,m}-\bs{t}^{l,m})\beta.
\end{equation}

\section{Examples} \label{sec:examples}
In this section, we evaluate and compare the PM finding algorithm with other DR methods on three different datasets:  a three-dimensional paraboloid, a swiss roll constructed with high additive noise (noisy swiss roll), and a simulation of collective behavior. Specifically, we use a standard version of eight different DR methods implemented in the Matlab Toolbox for DR \cite{van2007matlab} with default parameters listed in Table \ref{tab:comparison}. Note that we run each method in an unsupervised manner using the standard default values available in that toolbox.

\begin{table}[htp]
	\centering
	\caption{Default parameters of DR methods implemented in Matlab Toolbox for DR \cite{van2007matlab} for the paraboloid and noisy swiss roll examples.}
	\begin{tabular}{|l|l|}
		\hline
		Method & Parameters\\
		\hline
		 MDS & none\\ 
		 Isomap & nearest neighbors = 12\\
		 Diffusion maps & variance = 1 and operator on the graph = 1\\
		 KPCA & variance = 1 and Gaussian kernel is used\\
		 LSE & nearest neighbors = 12\\
		 LLE & nearest neighbors = 12\\
		 HLLE & nearest neighbors = 12\\
		 LE & nearest neighbors = 12 and variance = 1\\		
		\hline
	\end{tabular}
	\label{tab:comparison}
\end{table}

\subsection{Paraboloid}\label{sec:paraboloid}
\begin{figure*}[hpt]
	\includegraphics[width=1\linewidth]{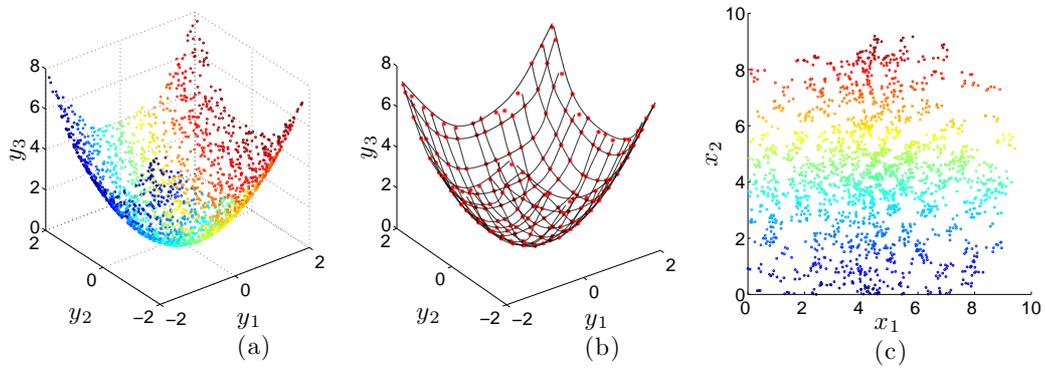}       	
        	\caption{Embedding data from a three-dimensional paraboloid with 2000 points (a) of raw data after smoothing (b) into two dimensional embedded space with intrinsic coordinate of the manifold (c). Points are colored to highlight the relative configuration.}
        	\label{fig:paraboloid_embedding}
\end{figure*}

\begin{figure*}[hpt]
	\includegraphics[width=1\linewidth]{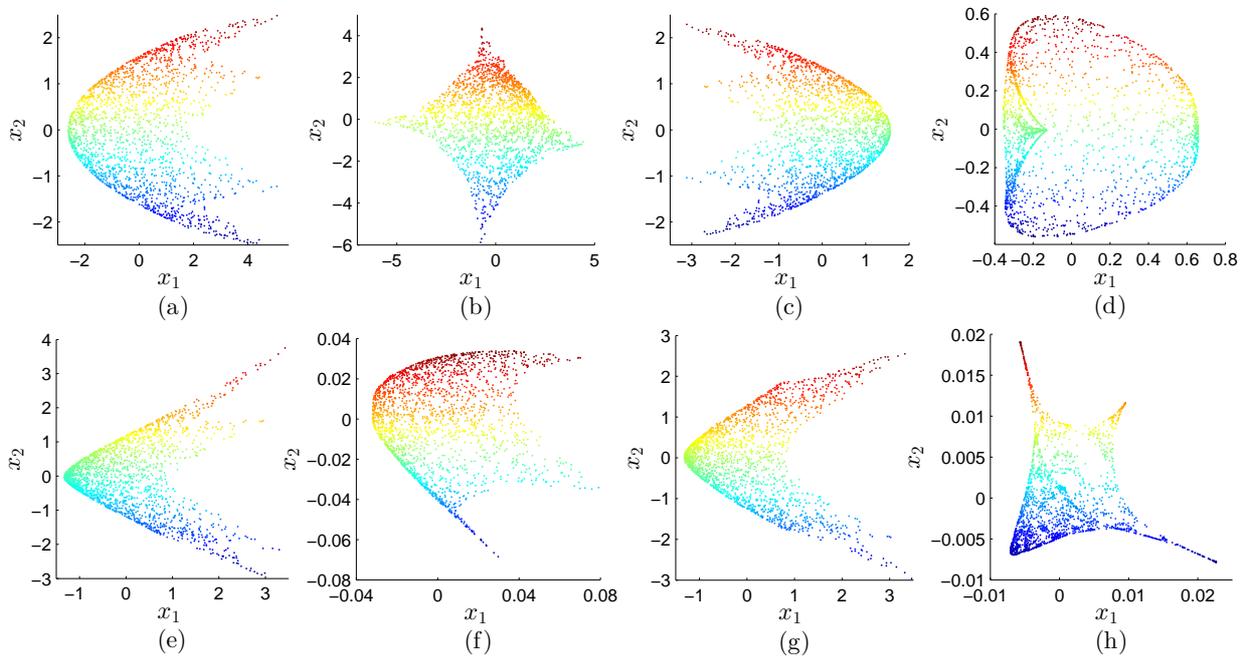}       	
	        	\caption{Two-dimensional embeddings of the paraboloid by MDS (a), Isomap (b), Diffusion maps (c), KPCA (d), LSE (e), LLE (f), HLLE (g), and LE (h)}
        	\label{fig:paraboloid_comparison}
\end{figure*}

The three-dimensional paraboloid (Fig. \ref{fig:paraboloid_embedding}a) is generated in the form of 2000 points using 
\begin{equation}
y_3=y_1^2+y_2^2+\epsilon,
\end{equation}
where $\bs{y}$ is a point in the three-dimensional space with $y_1, y_2 \in [-2, 2]$, and $\epsilon\sim\mathbb{U}(-0.05, 0.05)$ is a noise variable sampled from a uniform distribution ranging between -0.05 and 0.05. Since the data is generated along the $y_1,y_2$ axes, the two reference points are selected along the same just outside the range for each value. We run the clustering algorithm with the value of $n_{\mathcal{C}}=14$ for each reference point to generate corresponding sets of clusters with algorithm \ref{alg:cluster}. Next, we run the smoothing algorithm with a smoothing parameter $p=0.9$ to produce a set of representative splines for the data (Fig. \ref{fig:paraboloid_embedding}b). Finally, we embed the manifold by first finding the intersection points and, then, computing the distance of each point from the axis splines (Fig. \ref{fig:paraboloid_embedding}c). 

The two-dimensional embedding manifold generated using our algorithm preserves the topology of the data as shown in Figure \ref{fig:paraboloid_embedding}c. Furthermore, the axes in the embedding manifold retain the scale of the original data in terms of the distance between individual points. For a comparison, we run other NDR methods on this data set with parameters specified in table \ref{tab:comparison} (Fig.\ref{fig:paraboloid_comparison}). Results show that except Isomap and the proposed method, none of the NDR methods are able to preserve the two-dimensional square embedding after dimensionality reduction. Between Isomap and the PM approach, the PM approach retains the scale of the manifold as well.  
  
\subsection{A noisy swiss roll}
\begin{figure}[hpt]
	\centering
       	\includegraphics[width=1\textwidth]{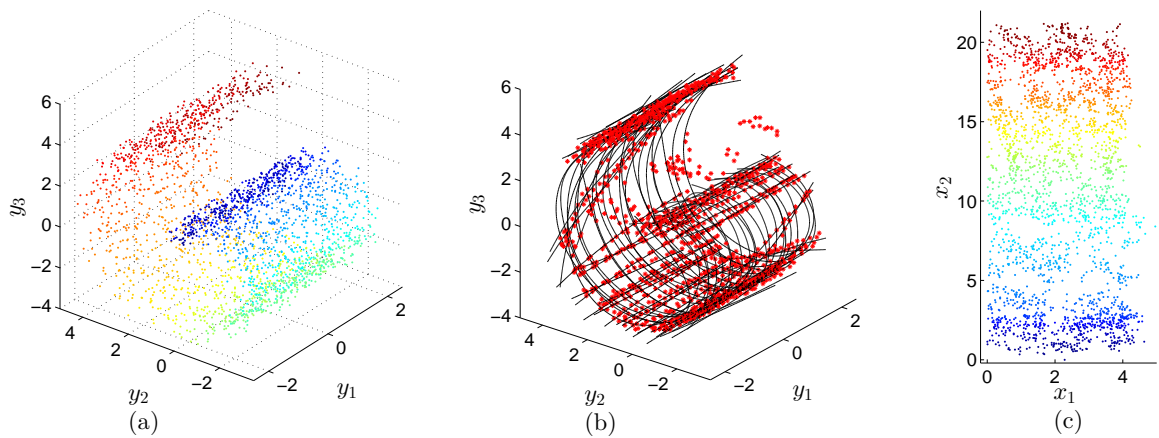}
	\caption{Two-dimensional embedding of noisy swiss-roll with 2500 points, raw data (a) of the noisy swiss roll is clustered and smoothened in order to obtain the principal manifold (b) and then embedded in two dimensions with intrinsic coordinates denoting the distance from the center of the manifold coordinate (c).}
	\label{fig:noised_swiss_roll_embedding}
\end{figure}

\begin{figure}[hpt]
	\centering
       	\includegraphics[width=1\textwidth]{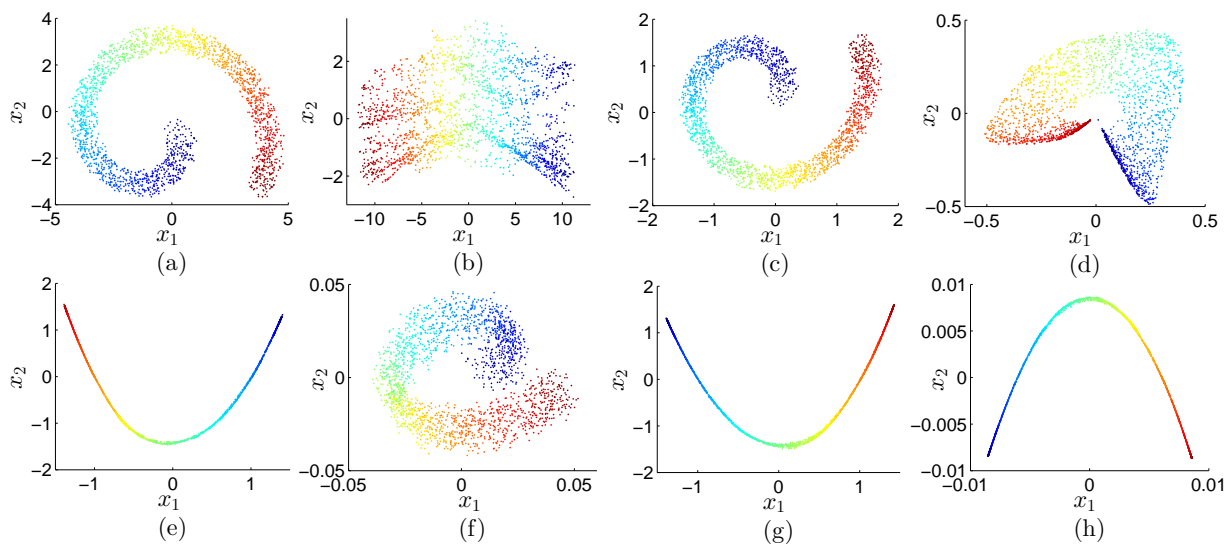}
	\caption{Two-dimensional embeddings of the noisy swiss-roll by MDS (a), Isomap (b), Diffusion maps (c), KPCA (d), LSE (e), LLE (f), HLLE (g), and LE (h).}
	\label{fig:noised_swiss_roll_comparison}
\end{figure}

In a second example, we run our algorithm on a 2500 point three-dimensional noisy swiss roll dataset with high additive noise (Fig. \ref{fig:noised_swiss_roll_embedding}a). The dataset is generated using 
\begin{equation}
\begin{split}
y_1=\theta^{0.8}\cos\theta+\epsilon \\
y_2=\theta^{0.8}\sin\theta+\epsilon \\
y_3=\phi+\epsilon
\end{split}
\label{eqn:noisy_sr}
\end{equation} 
where $\theta \in [0.25\pi, 2.5\pi], \ \phi \in [-2,2]$, and $\epsilon \sim\mathbb{U}(-0.4, 0.4)$  is a noise variable sampled from a uniform distribution ranging between -0.4 and 0.4. Two reference points are chosen along two major principal component directions $\bs{v}_1$ and $\bs{v}_2$ at a distance of $\sigma_1$ and $\sigma_2$ units away from the data centroid $\bs{\mu}$. We run the clustering algorithm with the value of $n_{\mathcal{C}}=15$ for each reference point, however, due to folds and therefore gaps in the data, the clusters along $\bs{v}_2$ are further split into 42 subclusters. The value of smoothing parameter $p=0.75$. Figure \ref{fig:noised_swiss_roll_embedding}b and \ref{fig:noised_swiss_roll_embedding}c show the smoothing splines in a grid pattern that are used to embed the manifold, and the resulting embedding. For a comparison, we run other NDR methods on this data set with parameters specified in table \ref{tab:comparison} (Fig. \ref{fig:noised_swiss_roll_comparison}).
Comparison between existing NDR methods show that while only Isomap and KPCA are able to flatten the swiss roll into two dimensions, Isomap preserves the general rectangular shape of the flattened swiss roll. A majority of NDR methods suppress one component out of the two principal components in the raw data revealing an embedding that resembles a one-dimesional curve. In contrast, the PM approach is able to embed the swiss roll into a rectangle while preserving the scale.

We note that by changing the value of the smoothing parameter we control for the level of noise-rejection in the original data, a feature that is not available in the existing DR algorithms \cite{van2009dimensionality}. Although the splines form illegal connections, our algorithm is still able to preserve the two-dimensional embedding; this is due to the inherent embedding step which is able to overcome such shortcuts while computing the low-dimensional coordinates. Furthermore, we note that our method is able to better exclude noisy data from the embedding  automatically in the form of outliers; neighborhood based algorithms would be mislead to attempt to include these as valid points \cite{Tenenbaum2000}. For example, the two-dimensional embedding obtained using Isomap shows an unforeseen bend demonstrating that, despite fewer outliers, the general  topological structure is compromised more heavily. This global anomaly is a direct consequence of the Isomap's attempt to include noisy points in the embedding; points that are otherwise ignored by the principal manifold.  

\subsection{Collective behavior: simulation of predator mobbing}
\begin{figure}
	\centering
	\includegraphics[width=.995\linewidth]{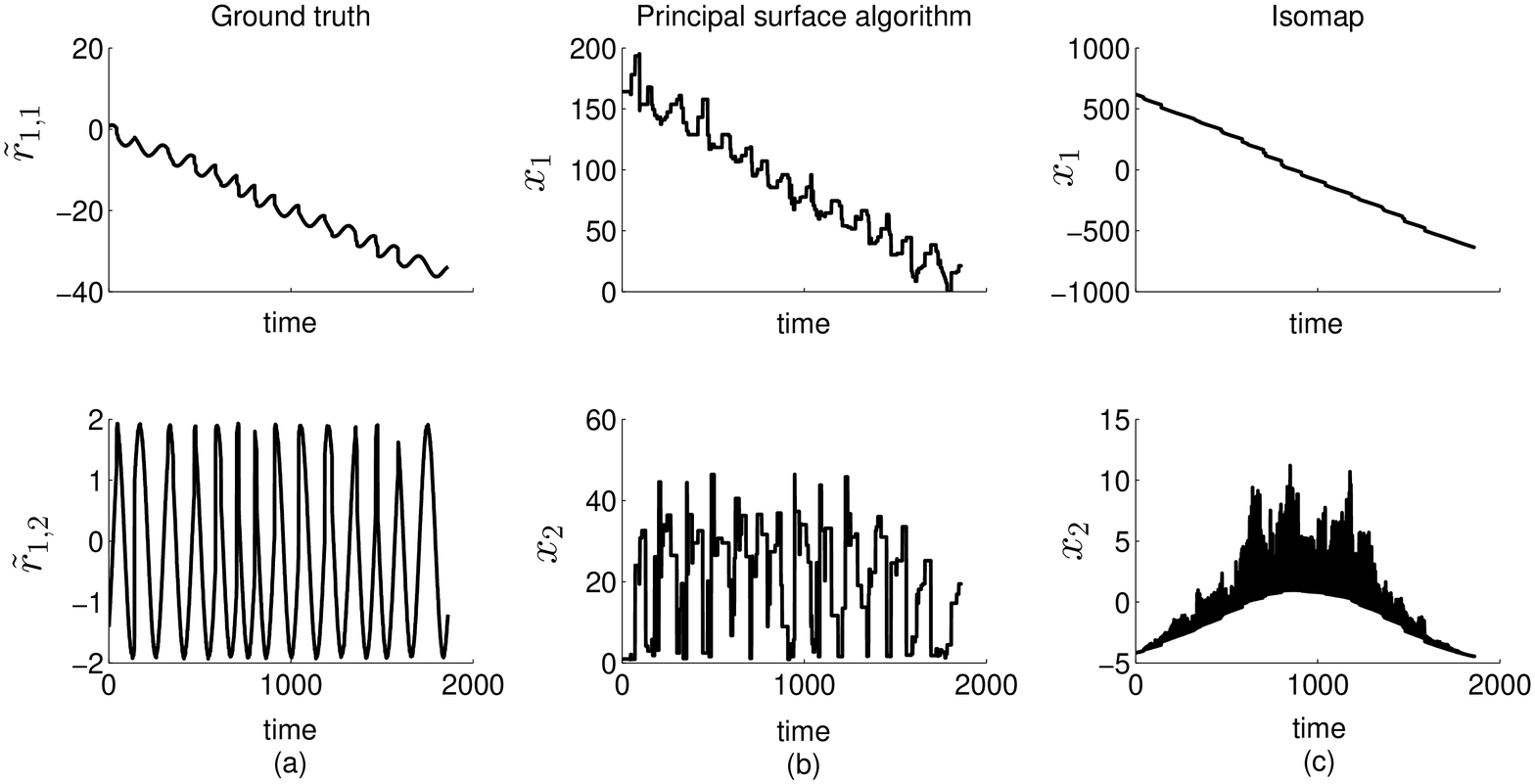}
	\caption{Embedding a forty dimensional dataset of an idealized simulation of predator mobbing where twenty particles move around a translating circle. A comparison of individual axes of the embedding produced by the principal surface algorithm described in this paper and those produced using Isomap shows that our algorithm performs distinctly better in terms of preserving the time signature as well as range of the axes. The nearest neighbor parameter used for Isomap is equal to 10.}
	\label{fig:predator_mobbing_algorithm2}
\end{figure}

In a third example of low-dimensional embedding, we generate a dataset representing a form of collective behavior. We simulate point-mass particles revolving on a translating circle to represent a form of behavior called predator induced mobbing \cite{Dominey1983}. In prey animals, predator mobbing, or crowding around a predator, is often ascribed to the purpose of highlighting the presence of a predator \cite{Dominey1983}. In our idealized model, we assume that $N$ agents ($N=20$) are revolving on the circumference of a half circle whose center marks the predator and translating on a two-dimensional plane which is orthogonal to the axis of the revolution. To generate a two-dimensional trajectory with $\rho$ revolutions around a translating center we represent the two-dimensional position of an agent  $i$ at $k$, $r_i[k]=\begin{bmatrix}y_{2i-1}[k], & y_{2i}[k]\end{bmatrix}^{\mathrm{T}}$, $\bs{r}_i[k] \in \mathbb{R}^2$
\begin{equation}
	\bs{r}_i[k]=s_i[k]\begin{bmatrix}\cos(2\pi \rho k/n_k +\pi i/N) \\
					\sin(2\pi \rho k/n_k + \pi i/N)
					\end{bmatrix} + k\bs{\mathrm{v}}_p[k]+ \bs{\epsilon}[k], 
\label{eqn:predator_mobbing}
\end{equation}
where $\bs{\mathrm{v}}_p[k]=[1, 1]^{\mathrm{T}}/80$ is the velocity of the predator $s_i[k] =3$ is the speed of the agent $i$, $n_k$ is the number of total time-steps, and $\bs{\epsilon}[k] \sim\mathbb{N}(\boldsymbol{0},\boldsymbol{1}/100)$ is the Gaussian noise variable with mean $\boldsymbol{0} \in \mathbb{R}^2$ and standard deviation $\boldsymbol{1}/100\in \mathbb{R}^2$. The first term of Eqn. (\ref{eqn:predator_mobbing}) assigns the relative positions of $N$ agents onto equally spaced points on the circumference of  a half circle of units $s_i[k]$, and the second term gives the velocity of the virtual center (predator) of the translating circle. The agent-dependent phase $\pi i/N$ creates spatial separation of individual members around the virtual center. Figure \ref{fig:predator_mobbing_isomap}a shows the resulting trajectories for 2000 time-steps with $\rho=14$ in a $d$-dimensional space (note that $d=2N$ in this case, where $N$ is the number of agents). The interaction between the agents  is captured in the form of noise $\bs{\epsilon}$. A high value means that the agents interact less with each other. 

To cluster the dataset, we use $n_{\mathcal{C}}=3$. The value of the smoothing parameter $p=0.9$. The embedding manifold is shown in Fig.\,\ref{fig:predator_mobbing_algorithm2}b with the start and end points of the trajectory. 
For comparison, we also run the Isomap algorithm on this dataset with a neighborhood parameter value set to 10. To ensure that the embedding is consistent and robust to our choice of neighborhood parameter, we run the algorithm with neighborhood parameters 5 and 15 and verify that the output is similar. Figure \ref{fig:predator_mobbing_algorithm2}c shows the resulting two-dimensional embedding coordinates each as a function of time.

Unlike the paraboloid and the noisy swiss roll, the predator mobbing dataset represents a dynamical system with a characteristic temporal evolution. Therefore, we use cross-correlation between the embedding of each method and a reference ground truth to compare the performance of our algorithm to Isomap (Fig. \ref{fig:predator_mobbing_isomap}). For ground truth, we transform the  original data by a clockwise rotation of $\pi/4$ such that 
\begin{equation}
\tilde{\bs{r}}_i[k]=\frac{1}{\sqrt{2}}\begin{bmatrix}1 & -1 \\ 1 & 1\end{bmatrix}\bs{r}_i[k]
\end{equation}
and use the position of the first agent $\tilde{\bs{r}}_1[k] =\begin{bmatrix}\tilde{r}_{1,1}[k], & \tilde{r}_{1,2}[k]\end{bmatrix}^{\mathrm{T}}$ (this value serves as a descriptive global observable of the predator mobbing agents separated by a constant phase, than, for example, the group centroid, where the revolving motion is suppressed due to the instantaneous two-dimensional arrangement). We then cross-correlate each component of the two-dimensional embedding signal $\begin{bmatrix}x_1, x_2\end{bmatrix}^{\mathrm{T}}$ with the corresponding component of the ground-truth and add the two values. We compute correlations of the embedded data of our algorithm and Isomap with the ground truth along the oscillating directions as $0.3987$ and $0.0074$. Thus, the correlation of our algorithm with the ground truth is fifty times  higher than the correlation of Isomap with the ground truth. Correlation $p$-values under our method and Isomap with the ground truth are $0.0000$ and $0.3888$ respectively, showing that under 95\% of confidence interval, the values from our method are related to the ones available from the ground truth. Figure \ref{fig:predator_mobbing_isomap} visually demonstrates the same where $x_2$ from the principal surface follows the same trend as $\tilde{r}_{1,2}$. In addition, the range of values of the embedding coordinates is similar to the original values. 

\section{Performance analysis}\label{sec:performance}
To analyze the performance of the algorithm, we use a distance-preserving metric between the original and the embedding data in terms of adjacency distance of graph connectivity. For both the original and the embedding data, we first search $k$-nearest neighbors through the algorithm given in \cite{friedman1977algorithm} and then produce two individual weighted graphs based on points connectivity in the data sets. A graph constructed through nearest neighbor search is a simple graph\footnote{A simple graph is an undirected graph that does not contains loops (edges connected at both ends to the same vertex ) and multiple edges (more than one edge between any two different vertices) \cite{balakrishnan2012textbook}.} that does not contain self-loops or multiple edges. We compute the adjacency distance matrix $A$ for the original data as
\begin{equation}
A_{i,j}  =
\begin{cases}
d(i,j) & : \exists \text{ an edge } ij \text{ in the graph of the original data}\\
0 & : \text{otherwise}
\end{cases}
\end{equation}
and the adjacency matrix $\tilde{A}$ for the embedding data as
\begin{equation}
\tilde{A}_{i,j}  =
\begin{cases}
d(i,j) & : \exists \text{ an edge } ij \text{ in the graph of the embedding data}\\
0 & : \text{otherwise}
\end{cases}
\end{equation}
Here, $A_{ij}$ and $\tilde{A}_{ij}$ are the $(i,j)-$th entries of the adjacency matrices $A$ and $\tilde{A}$, and $d(i,j)$ is the Euclidean distance between nodes $i$ and $j$. For $n$ points, this metric is computed as the normalized number of pairwise errors between entries of the adjacency distance matrices as
\begin{equation}\label{eqn:pairwise_error}
\Delta(\tilde{A},A)=\frac{1}{n k}\sum_{i,j}\mid\tilde{A}_{ij}-A_{ij}\mid,
\end{equation}
Equation (\ref{eqn:pairwise_error}) is a modification of the structure preserving metric in \cite{shaw2009structure} where the connectivity is replaced by the distance and the denominator $n^2$ with the number of edges $n k$. We study the dependence of $\Delta(\tilde{A},A)$ on the smoothing parameter, the noise in the dataset, and the data density.

\subsection{Error dependence on smoothing}
\begin{figure}[htp]
	\centering
	\includegraphics[width=.7\linewidth]{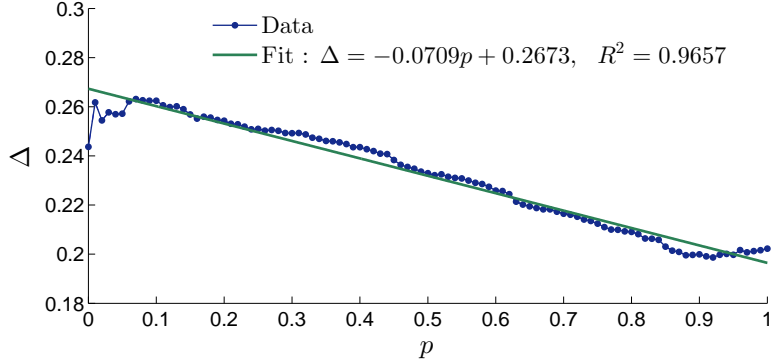}
	\caption{Variation of $\Delta$ with respect to $p$ and corresponding linear fit for a noise free swiss-roll. $\Delta$ is computed using Eq. (\ref{eqn:pairwise_error}) on raw and embedded data.}
	\label{fig:analysis_delta_p}
\end{figure}

The primary input to the algorithm described in this paper is the smoothing parameter that controls the degree of fitness and noise-rejection by the cubic smoothing splines. Figure \ref{fig:analysis_delta_p} shows the dependence of the distance preservation error $\Delta$ in term of adjacency distance on the smoothing parameter $p$. Using a noise-free swiss roll dataset comprising 3000 points in three dimensions, we vary the smoothing parameter between 0 and 1 with an increment of  0.01 while keeping the location of reference points and cluster width consistent. The error dependence shows a linear decay with $R^2$-value \footnote{In a fit of a data set, $R^2$-value (coefficient of determination) ranges between 0 and 1, and measures the goodness of the fit so that 1 is the best while 0 is the worst. For a set of points $\{y_i\}_{i=1}^n$ associated with fitted values $\{f_i\}_{i=1}^n$, coefficient of determination or $R^2$-value  is defined as, $1-\frac{\sum_i^n \left(y_i-f_i\right)^2}{\sum_i^n \left(y_i-\bar{y}\right)^2}$ where $\bar{y}=\sum_i^n y_i/n$.} of 0.9728 and becomes nearly constant at a $p=0.88$, beyond which the change in error for an increase in the value of $p$ is negligible. A higher value of $p$ improves the data-fit but a low value improves smoothness. Since the graph connectivity is dependent on the relative point configuration, it is expected that $\Delta$ will rise with increasing smoothness. On the other hand, the error dependence shows that a given degree of smoothing can be attained beyond the value of $p=0.88$ without an accompanied increase error. Thus, as a design parameter, the value $p=0.9$ may be used as a starting point for all datasets to investigate the embedding manifold.

\subsection{Error with noise}
\begin{figure}[htp]
	\centering
	\includegraphics[width=.7\linewidth]{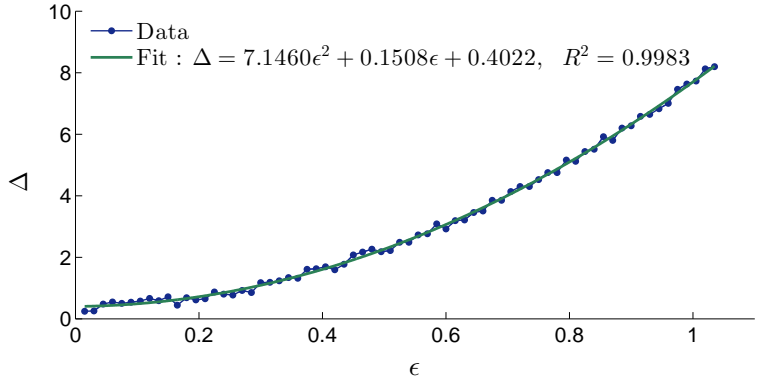}
	\caption{Variation of $\Delta$ with respect to noise ($\epsilon$) and corresponding quadratic fit  $\Delta$ is computed over adjacency matrices for the raw data and embedded data by Eq. (\ref{eqn:pairwise_error}).}
	\label{fig:analysis_delta_noise}
\end{figure}

To analyze the change in pairwise error with increase in noise, we create multiple same-sized swiss-roll datasets generated using (\ref{eqn:noisy_sr}) with noise value $\epsilon$ ranging between 0 and 1.035 with an increment of 0.015. The number of points in the dataset are 3000 and the value of smoothing parameter $p=0.9$. Figure \ref{fig:analysis_delta_noise} shows the variation of $\Delta$ with respect to $\epsilon$. The error $\Delta$ increases quadratically with $R^2$-value 0.9983 when noise increases. This relation shows that, the error of the embedding is affected quadratically by the intense of the associated noise in the data.

\subsection{Data density}
\begin{figure}[htp]
	\centering
	\includegraphics[width=.72\linewidth]{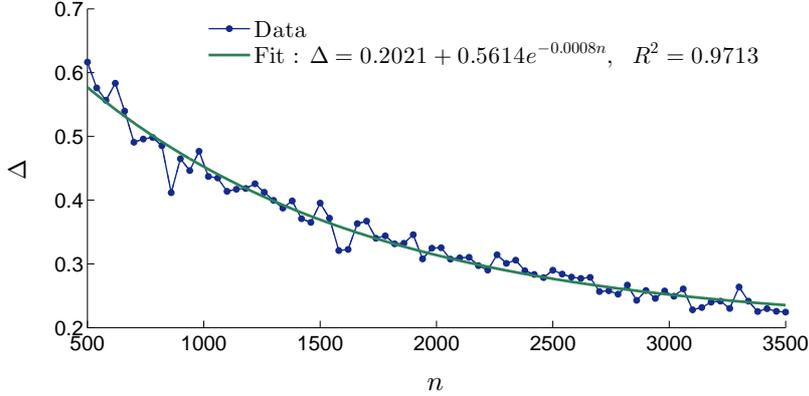}
	\caption{Variation of $\Delta$ with respect to number of points ($n$) and corresponding exponential fit with $p=0.9$ and $\epsilon=0.2$. $\Delta$ is computed over adjacency matrices for the raw data and embedded data by Eq.  (\ref{eqn:pairwise_error}). }
	\label{fig:analysis_delta_sparsity}
\end{figure}

To analyze the dependence of $\Delta$ on data density we sub-sample the swiss roll dataset such that the number of points are systematically decreased. The amount of noise is set at $\epsilon=0.2$ and a sequence of data sets are generated with number of points between 500 and 3500 with increment of 40. The value of the smoothing parameter $p$ is fixed at 0.9. We see that the pairwise error ($\Delta$) decreases exponentially with $R^2$-value $0.9713$ from an initial high value (Fig. \ref{fig:analysis_delta_sparsity}).

\section{Conclusion} \label{sec:conclusion}

In this paper, we introduce a PM finding dimensionality reduction algorithm that works directly on raw data by approximating the data in terms of cubic smoothing splines. We construct the splines on two sets of non-overlapping slices of the raw data created using parallel hyperplanes that are known distance apart. The resulting PM is a grid-like structure that is embedded on a two-dimensional manifold. The smoothness of the splines can be controlled via a smoothing parameter that selectively weighs fitting error versus the amount of noise rejection. The PM is represented by the intersection points between all pairs of splines, while the embedding coordinates are represented in the form of distance along these splines.

The spline representation improves regularity of an otherwise noisy dataset. We demonstrate the algorithm on three separate examples including a paraboloid, a noisy swiss roll, and a simulation of collective behavior. Upon embedding these datasets on two-dimensional manifolds, we find that in the case of the paraboloid, the algorithm is successful in preserving the topology and the range of the data. In the case of the noisy swiss roll, we find that the algorithm is able to reject high noise, thereby preserving the regularity in the original dataset. Unlike Isomap, our algorithm for finding PM is able to maintain a general structure. In other words, the algorithm fails gracefully giving the user enough indication of the trend of increasing error differently than Isomap, which fails abruptly and dramatically \cite{Mekuz2006}. This property has a distinct advantage over other methods that are more sensitive to input parameters and therefore make it difficult to identify when the algorithm is embedding correctly. In contrast, the proposed algorithm on the swiss-roll dataset shows that in its current form, is inefficient to embed the dataset with holes faithfully onto a two-dimensional space, since it is highly focused on revealing a smooth underlying manifold.

We used a simulation of twenty particles moving in a translating circle as a representation of a two-dimensional data embedded in a forty-dimensional state. In collective behavior, this form of motion approximates predator mobbing. Cross-correlation between the embedded signals and the ground-truth available in the form of a transformed signal shows that our algorithm is able to replicate the time-history on a low-dimensional space than for example Isomap. Furthermore, we see that the range of axes in the embedding manifold available from our algorithm is closer to that in the original data. 

The analysis of the distance preserving ability of our algorithm with respect to the smoothing parameter reveals that an increase in the smoothing parameter reduces the error in the structure. This is expected because a high value of $p$ implies that the cubic smoothing spline weights data fit more than smoothing. However, we also note that the amount of error saturates near $p=0.88$, showing that values close to but not equal to one may also be used to represent the data. We observe a similar trend with respect to the amount of noise introduced in a swiss-roll dataset, where we find that the amount of error increases quadratically  with an increase in noise. Finally, we find that the algorithm is able to work on sparse datasets up to a representation with five hundred points.

Operating and representing raw data for dimensionality reduction provides an opportunity for extracting true embeddings that preserve the structure as well as regularity of the dataset. By demanding a two-dimensional embedding we also provide a useful visualization tool to analyze high-dimensional datasets. Finally, we show that this approach is amenable to high-dimensional dynamical systems such as those available in dataset of collective behavior.

\section{Acknowledgements}
Kelum Gajamannage and Erik M. Bollt have been supported by the National Science Foundation under grant no. CMMI- 1129859. Sachit Butail and Maurizio Porfiri have been supported by the National Science Foundation under grants nos. CMMI- 1129820. 

\appendix
\section{Algorithms} \label{sec:algorithms}
Here, we state the algorithms of dimensionality reduction by principal manifold by three components as, clustering, smoothing, and embedding.
\begin{algorithm*}[htp]
\caption{ \textit{Clustering} : Data matrix $\mathcal{D}$ and number of clusters with respect to the first ($n_\mathcal{C}^{1}$) and second ($n_\mathcal{C}^{2}$) reference points  are three inputs in the clustering algorithm. As the output, this produces two sets of clusters $\mathcal{C}_j^i \ ; \ j=1, \ldots, n_\mathcal{C}^i$ for $i=1, 2$ from $\mathcal{D}$ such that one for each reference point.}\label{alg:cluster}
\begin{algorithmic}[1]
\Procedure{Clustering}{$\mathcal{D}, n_\mathcal{C}^{1}, n_{\mathcal{C}}^2$}
\State Compute the mean $\bs{\mu} \in \mathbb{R}^d$ of the input data $\mathcal{D}$.
\State Perform Principal Component Analysis (PCA) \cite{golub2012matrix, jolliffe2005principal} on the input data matrix $\mathcal{D}=\begin{bmatrix} \bs{y}_1| \bs{y}_2| \ldots |\bs{y}_n\end{bmatrix}$ to obtain two largest principal components $(\bs{v}_1, \sigma_1)$ and $(\bs{v}_2, \sigma_2)$, where $\bs{v}_i$ is the $d$-dimensional coefficients and $\sigma_i$ is the eigenvalue of the $i$-th PC for $i=1, 2$
\State In order to assure that two reference points are in the directions of two PCs from the mean, compute the first reference point $\bs{q}_1=\bs{\mu}+\bs{v}_1\sigma_1$, and the second $\bs{q}_2=\bs{\mu}+\bs{v}_2\sigma_2$.
\For{each reference point $\bs{q}_i ; i=1,2$}
\State Segment the line joining the reference points $\bs{q}_i=\bs{\mu}+\bs{v}_i\sigma_i$ and the point $\tilde{\bs{q}}_i=\bs{\mu}-\bs{v}_i\sigma_i$ into $n_\mathcal{C}^{i}$ parts with equal width using the ratio formula given in equation (\ref{eqn:line_partition}) to obtain a set of points $\bs{a}_j; j=0, \ldots, n_\mathcal{C}^i$  \cite{Protter1988}, where $\bs{a}_0=\tilde{\bs{q}}_i$ and $\bs{a}_{n_\mathcal{C}^i}=\bs{q}_i$
\State For $j=1, \ldots, n_\mathcal{C}^i$, choose data between hyperplanes, which are made through points $\bs{a}_{j-1}$ and $\bs{a}_{j}$ normal to the line joining $\tilde{\bs{q}}_i$ and $\bs{q}_i$, into the cluster $\mathcal{C}^i_j$ by satisfying the inequalities in (\ref{eqn:enq_clustering}).
\EndFor
\EndProcedure
\end{algorithmic}
\end{algorithm*}

\begin{algorithm*}
\caption{ \textit{Smoothing} : Smoothing algorithm produces two sets of cubic smoothing splines with respect to both reference points to represents the data in clusters. This algorithm inputs spline smoothing parameter $p$ and two sets of clusters $\mathcal{C}_j^i \ ; \ j=1, \ldots, n_\mathcal{C}^i$ for $i=1, 2$ made in the clustering algorithm, and outputs two sets of cubic  smoothing splines  $\bs{S}_j^i \ ; \ j=1, \ldots, n_\mathcal{C}^i$ for $i=1, 2$.}\label{alg:smooth}
\begin{algorithmic}[1]
\Procedure{Smoothing}{$p$; $\mathcal{C}^i_j, j=1, \ldots, n_\mathcal{C}^i$ for $i=1, 2$}
\For{each reference point $\bs{q}_i ; i=1,2$}
\For{all clusters $\mathcal{C}^i_j;j=1, \ldots, n_\mathcal{C}^i$}
\State Compute the neighborhood graph using range-search with the distance set as the cluster width $2\sigma_i/n_{\mathcal{C}}^i$.
\State Compute the longest geodesic $\mathcal{G}^i_j$ using Dijkstra's algorithm \cite{dijkstra1959note, leiserson2001introduction}. $\mathcal{G}^i_j$ is a set of points, where each point is in $\mathbb{R}^d$.
\State For points in $\mathcal{G}^i_j$, use (\ref{eqnCSS}) and the value of the smoothing parameter $p$ to produce a smoothing spline $\bs{S}^i_j \in \mathbb{R}^d$ representation for data in the cluster $\mathcal{C}^i_j$.
\EndFor
\EndFor
\EndProcedure
\end{algorithmic}
\end{algorithm*}

\begin{algorithm*}
\caption{ \textit{Embedding} : This produces a grid structure by approximating the intersections of splines, followed by doing the embedding of a new point based on this structure. Embedding algorithm inputs a new point $z \in \mathbb{R}^d$ to embed and two sets of cubic  smoothing splines $\bs{S}_j^i \ ; \ j=1, \ldots, n_\mathcal{C}^i$ for $i=1, 2$ produced by smoothing algorithm, and this outputs two dimensional embedding coordinates $\left[x_1, x_2\right]^{\mathrm{T}}$ of $\bs{z}$.}\label{alg:embed}
\begin{algorithmic}[1]
\Procedure{Embedding}{$\bs{S}_j^i, j=1, \ldots, n_\mathcal{C}^i$ for $i=1, 2$ ; a point $\bs{z} \in \mathbb{R}^d$ for embedding}
\For{all pairs $(l,m)$ of smoothing splines $\bs{S}^1_l \times \bs{S}^2_m$ belonging to reference points 1 and 2}
\State Approximate the minimum distance between the two splines after discretizing each spline.
\State Choose the midpoint between the two closest points of the splines as the intersection point $\bs{t}^{l,m} \in \mathbb{R}^d$.
\EndFor
\State Pick a random intersection point as the origin $O \in \mathbb{R}^d$, and smoothing splines corresponding to the origin as axis splines.
\State Find the closest intersection point $\tilde{\bs{z}}$ for $\bs{z}$ in terms of Euclidean distance by equation (\ref{eqn:closest_intersection}), and the tangents to the splines at this point are called the local spline directions $(\bs{u}^1_{\bs{z}}, \bs{u}^2_{\bs{z}})$.
\State Compute the distances $\left[\tilde{x}_1, \tilde{x}_2\right]^{\mathrm{T}}$ from the manifold origin $O$ to $\tilde{\bs{z}}$ along axis splines by equation (\ref{eqn:axis_distance}).
\State Project the vector $\tilde{\bs{z}}-\bs{z}$ onto the local coordinate system created using the spline directions as equation (\ref{eqn:local_distance}) at the intersection point $\tilde{\bs{z}}$ and find the projections $\left[\delta x_1, \delta x_1\right]^{\mathrm{T}}$. The final embedding coordinates are given as $x_1=\tilde{x}_1+ \delta x_1$, $x_2=\tilde{x}_2+ \delta x_2$.
\EndProcedure
\end{algorithmic}
\end{algorithm*}

\newpage
\section{Computational complexity}\label{sec:complexity}
The three steps of the algorithm, namely, clustering, smoothing, and embedding have different computational complexities. In particular, for $n$ points with dimensionality $d$, the clustering algorithm has a complexity 
\begin{equation}
O\left(m^3\right) ; \ \text{ where } \ m=\mathrm{min}(n, d),
\label{eqn:clustering_complexity}
\end{equation}
which is dominated by the PCA \cite{ jolliffe2005principal, jackson2005user} step.  

The computational complexity of the smoothing algorithm is dominated by the  Dijkstra's algorithm \cite{dijkstra1959note} for calculating the longest geodesics. Here, we first assume that each cluster with respect to the first reference point contains an equal number of points. Thus, partitioning a data set of $n$ point into $n_{\mathcal{C}}^1$ clusters yields $n/n_{\mathcal{C}}^1$ points in a cluster. The complexity of generating the longest geodesic in each such cluster is $O\left(n/n_{\mathcal{C}}^1\log \left(n/n_{\mathcal{C}}^1\right)\right)$ \cite{dijkstra1959note, leiserson2001introduction}, which sums-up $n_{\mathcal{C}}^1$ times and implies the total complexity $O\left(n\log \left(n/n_{\mathcal{C}}^1\right)\right)$ of making all geodesics with respect to the first reference point. Similarly, for the second reference point, the same procedure is followed for all $n_{\mathcal{C}}^2$ clusters, each containing $n/n_{\mathcal{C}}^2$ points, to obtain a complexity of $O\left(n\log \left(n/n_{\mathcal{C}}^2\right)\right)$. Altogether, geodesics from both reference points contributes a computational complexity of
\begin{equation}
O\left(n\log \left(\frac{n^2}{n_{\mathcal{C}}^1 n_{\mathcal{C}}^2}\right)\right)
\end{equation}
for the smoothing algorithm. 

In the embedding algorithm, discretizing splines and approximating the minimum distance are dominant in the total cost of that algorithm. Here, we first calculate the mean length of splines with respect to the first reference point, denoted as $\tilde{\bs{S}}^{1}\in \mathbb{R}^d$, and the one for the second reference point, denoted as $\tilde{\bs{S}}^{2}\in \mathbb{R}^d$. Then, the computational cost associated with approximating intersection of these two splines are computed. For a pair of $d$-dimensional splines, computation of the Euclidean distance between any two points, one from each spline, has complexity of $3d$ from the operations of subtraction, squaring, and summation along all $d$ dimensions. Since each spline $\tilde{\bs{S}}^{1, 2}$ is discretized at a mesh size $h$, while $\tilde{\bs{S}}^{1}$ has $\tilde{\bs{S}}^{1}/h$ points on it, $\tilde{\bs{S}}^{2}$ has $\tilde{S}^{2}/h$ points. The computational cost of approximating closest distance between any two points one from each spline $\tilde{\bs{S}}^{1, 2}$ is $3d\tilde{\bs{S}}^1 \tilde{S}^2/h^2$. Here, for simplicity, we assume that the length of each spline with respect to the first and second reference points have same lengths as $\tilde{\bs{S}}^{1}$ and $\tilde{S}^{2}$ consecutively. Thus, the total complexity of the embedding algorithm, which has $n^1_\mathcal{C}$ and $n^2_\mathcal{C}$  number of splines with respect to the first and second reference points, is 
\begin{equation}
O\left(\frac{3d}{h^2} n^1_\mathcal{C} n^2_\mathcal{C} \tilde{\bs{S}}^1 \tilde{\bs{S}}^2\right)
\label{eqn:embedding_complexity}
\end{equation}
Clustering, embedding algorithms and steps $1-6$ in the embedding algorithm are performed only once for a given data set, thus after they are completed, $n$ new points are embedded with $O(n)$.





\bibliographystyle{elsarticle-num}







\end{document}